\documentclass[russian]{report}
\usepackage[T2A]{fontenc}
\usepackage[cp1251]{inputenc}
\usepackage[russian]{babel}
\usepackage{amsmath}
\usepackage{mathtext}
\usepackage{euscript}
\usepackage{amsfonts}
\usepackage{amssymb}
\usepackage{amscd}
\begin{document}

\begin{center}
{\bfseries\large The geometry of a deformation of the standard addition on the integral lattice.}
\end{center}

\begin{center}
Stanislav Y. Tsarev
\end{center}
\par\medskip
\par\medskip

{\footnotesize	{\bfseries Abstract.} Let $\mathfrak A_n$ be the subset of the standard integer lattice $\mathbb Z^n$, $\mathfrak A_n\subset\mathbb Z^n$ which is defined by the condition $\mathfrak A_n=((a_1,\ldots,a_n)\in\mathbb Z^n |\ a_i\not\equiv a_j\mod n,\ \forall i,j\in \{1,\ldots\,n\})$. It is clear that the standard addition on the lattice $\mathbb Z^n$ doesn’t induce the group structure on the set $\mathfrak A_n$ since the componentwise sum of some two vectors may contain components which are equal modulo $n$. Our aim is to find a new associative multiplication on the lattice $\mathbb Z^n$ such that the induced multiplication on the set $\mathfrak A_n$ gives it the group structure. In this paper the group structure on the subset $\mathfrak A_n$ of the integer lattice $\mathbb Z^n$ is studied by means of the constructions, which are described in [1]. The geometric realization of this group in the enveloping space and its generators and relations between them are found. We begin with the main constructions and the results we need for them.
}
\par\medskip
\par\medskip

{\bfseries\Large 1. Introduction and the main definitions}
\par\medskip
\par\medskip

	Let $G$ be a group with multiplication $m(g_2,g_1)=g_2g_1$. For a space $V$ with a right action $\alpha$ of the group $G$, let us introduce the multiplication $m_\alpha$ in the space of maps $G^V$:
$$m_\alpha:\ \phi_2\ast\phi_1(v)=\phi_2(v)\phi_1(v\phi_2(v)),$$
where $\phi_1,\ \phi_2\in G^V$, $vg=\alpha(v,g)$. 
\par\medskip

	Obviously, in the case when $\alpha$ acts trivially, the corresponding multiplication coincides with the pointwise multiplication in $G^V$.
\par\medskip

	Denote by $G^V_\alpha$ the space with multiplication $m_\alpha$. Let $i\colon G\to G^V$ be the map that takes the point g to the constant map to this point.
\par\medskip

	In what follows we need the lemma from [1], the proof of which is provided for the convenience of the reader.
\par\medskip 

	{\bfseries Lemma 1.} 1) The set $G^V_\alpha$ is the semigroup with identity with respect to the multiplication, which was introduced above
\par\medskip

2) The map $i$ is the homomorphism from $G$ to $G^V_\alpha$.
\par\medskip

{\bfseries Proof.} Let us check associativity for the multiplication $m_\alpha$:
$$(\phi_3\ast(\phi_2\ast\phi_1))(v)=\phi_3(v)(\phi_2\ast\phi_1)(v_1)=\phi_3(v)(\phi_2(v_1)\phi_1(v_2)),$$
where $v_1=v\phi_3(v)$ and $v_2=v_1\phi_2(v_1)$. On the other hand,
$$((\phi_3\ast\phi_2)\ast\phi_1)(v)=(\phi_3\ast\phi_2)(v)\phi_1(v(\phi_3\ast\phi_2(v)))=(\phi_3(v)\phi_2(v_1))\phi_1(v_2).$$
The associativity of the multiplication in $G$ implies the equality of two resulting expressions. The identity in $G^V_\alpha$ is the element $id=i(e)$. Indeed, let us find:
$$id\ast\phi(v)=id(v)\phi(v\cdot id(v))=e\phi(ve)=\phi(v).$$ 
Similarly, in the reverse sequence
$$\phi\ast id(v)=\phi(v)id(v\phi(v))=\phi(v)e=\phi(v),$$
in what follows that the expressions are the same. The first assertion of the lemma is proved.
\par\medskip

	To prove the second assertion we check that $i$ is the homomorphism from $G$ to $G^V_\alpha$: $i(g_2)\ast i(g_1)(v)=i(g_2)(v)i(g_1)(v\cdot i(g_2)(v))=g_2g_1$. {\bfseries The lemma is proved.}
\par\medskip

	{\bfseries Lemma 2.} Let $V$ be a finite set. Then $\phi\colon V\to G$ is invertible in $G^V$ if and only if $v\phi(v)\colon V\to V$ is the bijection.
\par\medskip

{\bfseries Proof.} A mapping $\phi^{-1}\in G^V$ is the inverse of the $\phi\in G^V$, if $\phi\ast\phi^{-1}=\phi^{-1}\ast\phi=id$. According to the definition, $\phi\ast\phi^{-1}(v)=\phi(v)\phi^{-1}(v\phi(v))=e$. If $v\phi(v)\colon V\to V$ is the bijection, then for any $v\in V$ there is a unique $v^{\prime}$, such that $v=v^\prime\phi(v^\prime)$. So $\phi^{-1}$ is defined uniquely $\phi^{-1}(v)=(\phi(v^\prime))^{-1}$ and $\phi^{-1}$ is the right inverse. It is trivial to verify, that $\phi^{-1}$ is also the left inverse. Clear, that if $v\phi(v)\colon V\to V$ is not bijection, then $\phi^{-1}$ is not defined uniquely. {\bfseries The lemma is proved.}
\par\medskip
\par\medskip

{\bfseries\Large 2. The group description. Generators and retalions.}
\par\medskip
\par\medskip

	The main problem is to find the subgroup of invertible elements in the semigroup $G^V$ in a particular case. Namely, as a group $G$ we take the group $\mathbb Z$ of integer numbers under addition, and $V$ be the finite set, which consist of $n$ elements. An action $\mathbb Z$ on $V$ is given by a single permutation $\tau$, which corresponds to the unit $1\in\mathbb Z$; then number two corresponds to the permutation $\tau^2$, number three to $\tau^3$ and so on. For the permutation $\tau$ we take cyclic permutation $\tau=\begin{pmatrix}
1 & 2 & 3 & \ldots & n\\
n & 1 & 2 & \ldots & n-1
\end{pmatrix}
$, that is $\tau(1)=n,\ \tau(2)=1,\ \tau(3)=2,\ \ldots,\tau(n)=n-1$. For the convenience of writing any permutation we omit the top line and write only the lower.
\par\medskip

	Let us find out how the multiplication on $\mathbb Z^V$ looks like. Any function $\phi\colon V\to\mathbb Z$ is given by $n$ integer numbers; let’s multiply two functions by the rule which was described above: $(a_1\ldots ,a_n)\ast(b_1,\ldots b_n)=(a_1+b_{1\cdot\tau^{a_1}},a_2+b_{2\cdot\tau^{a_2}},\ldots ,a_n+b_{n\cdot\tau^{a_n}})=(a_1+b_{1+\widehat{-a_1}},a_2+b_{1+\widehat{1-a_2}},\ldots ,a_n+b_{1+\widehat{n-a_n-1}})$, where $\hat a=a\mod n$, i.e. the remainder of the division of $a$ by $n$.
\par\medskip

	Denote by $\mathfrak A_n$ the subgroup of invertible elements in $\mathbb Z^V$. 
\par\medskip

{\bfseries Assertion 1.} $(a_1,\ldots,a_n)\in\mathfrak A_n$ if and only if the numbers $1+\widehat{-a_1},\ldots,1+\widehat{n-1-a_n}$ are destinct.
\par\medskip

{\bfseries Proof.} According to lemma $2$, the numbers $\{1\phi(1),\ldots,n\phi(n)\}=\{1\cdot\tau^{a_1},\ldots,n\cdot\tau^{a_n}\}=\{1+\widehat{-a_1},\ldots,1+\widehat{n-1-a_n}\}$  must be distinct. {\bfseries The assertion is proved.}
\par\medskip

{\bfseries Assertion 2.} The group $\mathfrak A_n$ is isomorphic to the semidirect product of $\mathbb Z^n$ and symmetric group $S_n$, i.e. $\mathfrak A_n\cong\mathbb Z^n\rtimes S_n$, where the group $S_n$ acts on $\mathbb Z^n$ by permutations.
\par\medskip

{\bfseries Proof.} Let $(a_1,\ldots,a_n)\in\mathfrak A_n$, what is equivalent to the numbers $1+\widehat{-a_1},\ldots,1+\widehat{n-1-a_n}$ are distinct. It is clear, that the number of all such vectors up to the modulo $n$ is $n!$. Translating all these vectors by $(0,n-1,\ldots,2,1)$, we obtain $(a_1,a_2+n-1,\ldots,a_{n-1}+2,a_n+1)$. It is easy to check, that the coordinates of the obtained vectors are distinct modulo $n$. Indeed, the difference between any two components is $a_i-i-a_j+j$; on the other hand, the numbers $1+\widehat{i-1-a_i}$ and $1+\widehat{j-1-a_j}$ are distinct modulo $n$, therefore $i-1-a_i-j+1+a_j$ modulo $n$ is $a_i-i-a_j+j$ modulo $n$ and doesn’t equal zero. The new vectors represent all the points in $\mathbb Z^n$, the coordinates of which modulo $n$ are the numbers $1,2,\ldots,n-1,0$ in some order.
\par\medskip

	Now find out the structure of multiplication on the new vectors. The structure of multiplication must be maintained, that is the product of any two vector, moved by the displacement vector, is equal to the product of two vectors, which are the result of the translation by the displacement vector. Therefore, for the new vectors the following equality holds  $(a_1,\ldots,a_n)\times(b_1\ldots,b_n)=(a_1,a_2-n+1\ldots,a_n-1)\ast(b_1,b_2-n+1,\ldots,b_n-1)+(0,n-1,\ldots,1)=(a_1-\widehat{a_1}+b_{1+\widehat{-a_1}},a_2-\widehat{a_2}+b_{1+\widehat{-a_2}},\ldots,a_n-\widehat{a_n}+b_{1+\widehat{-a_n}})$.
\par\medskip

	Now we construct the homomorphism $\phi\colon\mathfrak A_n\to\mathbb Z^n\rtimes S_n$. Every vector in $\mathfrak A_n$ can be uniquely represented in the form $(nm_1+l_1,\ldots,nm_n+l_n)$. Let the map $\phi$ takes this vector to the pair $(z,s)\in\mathbb Z^n\rtimes S_n$, where $z=(m_1,\ldots,m_n)$, $s=(1+\widehat{-l_1},\ldots,1+\widehat{-l_n})$. We show that such correspondence is the homomorphism. Let’s multiply two vectors in $\mathfrak A_n$, $(nm_1+l_1,\ldots,nm_n+l_n)\times(nk_1+t_1,\ldots,nk_n+t_n)=(nm_1+nk_{1+\widehat{-l_1}}+t_{1+\widehat{-l_1}},\ldots,nm_n+nk_{1+\widehat{-l_n}}+t_{1+\widehat{-l_n}})$. The map $\phi$ takes this vector to the pair $((m_1+k_{1+\widehat{-l_1}},\ldots,m_n+k_{1+\widehat{-l_n}}),(1+\widehat{-t_{1+\widehat{-l_1}}},\ldots,1+\widehat{-t_{1+\widehat{-l_n}}}))$. On the other hand, the map $\phi$ takes the initial vectors to the pairs $((m_1,\ldots,m_n),(1+\widehat{-l_1},\ldots,1+\widehat{-l_n}))$ and $((k_1,\ldots,k_n),(1+\widehat{-t_1},\ldots,1+\widehat{-t_n}))$. Their product in the group $\mathbb Z^n\rtimes S_n$ is $((m_1+k_{1+\widehat{-l_1}},\ldots,m_n+k_{1+\widehat{-l_n}}),(1+\widehat{-t_{1+\widehat{-l_1}}},\ldots,1+\widehat{-t_{1+\widehat{-l_n}}}))$, therefore $\phi$ is the homomorphism. The map $\phi$ is obviously monomorphic and epimorphic, so $\phi$ is isomorphism. {\bfseries The assertion is proved.}
\par\medskip

	Let us now find the generators and relation in the group $\mathfrak A_n\cong\mathbb Z^n\rtimes S_n$. The idea behind the prood is that this group is easily described by three generators, and then the description is reduced to two generators.
\par\medskip

	{\bfseries Assertion 3.} For $n\ge 4$ the group $\mathbb Z^n\rtimes S_n$ may be described as follows 

$$
\left\{
\begin{array}{rcl}
\sigma^2&=&e\\
(\sigma\tau\sigma\tau^{-1})^3&=&e\\
(\sigma\tau^m\sigma\tau^{-m})^2&=&e,\ 2\le m\le n-2\\
(\sigma\tau)^{n-1}&=&\tau^n\\
\gamma\tau^k\sigma\tau^{-k}&=&\tau^k\sigma\tau^{-k}\gamma,\ 0 \le k \le n-3\\
\gamma\tau^l\gamma\tau^{-l}&=&\tau^l\gamma\tau^{-l}\gamma, \ 1 \le l \le n-1\\
\end{array}
\right.
$$ 
\par\medskip

	At first, we need the following lemma:
\par\medskip

	{\bfseries Lemma 3.} The symmetric group $S_n$ may be described as follows 
$$\{\tau,\sigma|\sigma^2,(\sigma\tau\sigma\tau^{-1})^3,(\sigma\tau^m\sigma\tau^{-m})^2,2 \le m \le n-2,(\sigma\tau)^{n-1}=\tau^n\}.$$
\par\medskip
	
	{\bfseries Proof.} It is known that the symmetric group $S_n$ has the next presentation
$$\{\sigma_1,\sigma_2,\ldots,\sigma_{n-1}|\sigma_i^2,(\sigma_i\sigma_{i+1})^3,\sigma_i\sigma_j=\sigma_j,\sigma_i,|i-j|>1\},$$
where $\sigma_i$ transposes $i$ and $i+1$.
\par\medskip

	As the new generators we take the following two permutations:
$$
\left\{
\begin{array}{rcl}
\sigma &=&\sigma_1\\
\tau &=&\sigma_1\sigma_2\ldots\sigma_{n-1}.\\
\end{array}
\right
.$$
\par\medskip

	It is easy to see that $\tau$ is the cycle $\begin{pmatrix}
1 & 2 & 3 & \ldots & n\\
n & 1 & 2 & \ldots & n-1
\end{pmatrix}
$. Let us deduce the complete system of relations for the new generators.
\par\medskip

	First we need to express the old generators via the new generators. We prove that $\sigma_i=\tau^{i-1}\sigma\tau^{1-i}$. This follows from the equality $\\ (\sigma_1\ldots\sigma_{n-1})\sigma_i(\sigma_{n-1}\ldots\sigma_1)=\sigma_1\ldots\sigma_i\sigma_{i+1}\sigma_i\sigma_{i+1}\sigma_i\ldots\sigma_1=\\ \sigma_1\ldots\sigma_{i-1}\sigma_{i+1}\sigma_i\sigma_{i+1}\sigma_{i+1}\sigma_i\sigma_{i-1}\ldots\sigma_1=\sigma_1\ldots\sigma_{i-1}\sigma_{i+1}\sigma_{i-1}\ldots\sigma_1=\sigma_{i+1}.$
\par\medskip

	Rewrite now the group presentation by the new generators. The relation $\sigma_i^2=e$ turns to $\sigma^2=e$, where $e$ is the unit. The relation $(\sigma_i\sigma_{i+1})^3=e$ turns to $(\sigma\tau\sigma\tau^{-1})^3=e$, and $\sigma_i\sigma_j=\sigma_j\sigma_i$ turns to $(\sigma\tau^m\sigma\tau^{-m})^2=e$, $2 \le m \le n-2$. The last relation is obtained from the formlae that express one generating system through another $(\sigma\tau)^{n-1}=\tau^n$. So we obtain the following presentation of the group $S_n$
$$\{\tau,\sigma|\sigma^2,(\sigma\tau\sigma\tau^{-1})^3,(\sigma\tau^m\sigma\tau^{-m})^2,2 \le m \le n-2,(\sigma\tau)^{n-1}=\tau^n\}.$$ {\bfseries The lemma is proved.}
\par\medskip

	{\bfseries Notation 1.} From the obtained relations for the system of generators for the group $S_n$ it is easy to deduce that $\tau^n=e$.
\par\medskip 

	Let us now find the generators and relation of the group $\mathbb Z^n\rtimes S_n$. As a generators one can take three elements $\sigma=(0,\ldots,0),(213\ldots n)$, $\tau=(0,\ldots,0),(n12\ldots n-1)$ and $\gamma=(0,\ldots,1),(123\ldots n)$. It is trivial to verify that they generate the whole group.
\par\medskip

	The relations between the elements $\sigma$ and $\tau$ were found in the previous lemma. It is easy to find another tow relations, namely $\gamma\tau^k\sigma\tau^{-k}=\tau^k\sigma\tau^{-k}\gamma$, $0 \le k \le n-3$ and $\gamma\tau^l\gamma\tau^{-l}=\tau^l\gamma\tau^{-l}\gamma$, $1 \le l \le n-1$. We show that the resulting system of the relations is complete.
\par\medskip

	{\bfseries Notation 2.} One can generalize two obtained relations $\gamma^a\tau^l\gamma^b\tau^{-l}=\tau^l\gamma^b\tau^{-l}\gamma^a$, $\gamma^a\tau^k\sigma\tau^{-k}=\tau^k\sigma\tau^{-k}\gamma^a$.
\par\medskip

	Let there is any relation of the form $\gamma^{\alpha}\tau^{\beta}\sigma\ldots\gamma^{\delta}\ldots=e$. If in the given relation $\tau$ is absent, i.e. the equality takes the form $\gamma^{\alpha}\sigma\gamma^{\beta}\sigma\ldots=e$, then by $\gamma\sigma=\sigma\gamma$ we obtain $\gamma^{\alpha+\beta+\ldots}\sigma^t=e$. This equality is true if and only if $t$ is even and $\alpha+\beta+\ldots=0$. Obviously, in this case the obtained relation is the consequence of the relation $\gamma\sigma=\sigma\gamma$.
\par\medskip

	Let the initial relation doesn’t contain $\sigma$, that is it takes the form $\gamma^{a_1}\tau^{b_1}\ldots\gamma^{a_p}\tau^{b_p}=e$. The numbers $a_1,\ldots,a_p$ are nonzero, and $b_1,\ldots,b_p$ are not divisible by $n$, otherwise one can reduce the length $p$.
\par\medskip

	Note that $p\ge 2$. Indeed, $\gamma^{a_1}\tau^{b_1}=(0,\ldots,a_1),\tau^{b_1}$, the first component is nonzero, so the initial equality isn’t true for $p=1$. 
\par\medskip

	For $p=2$ the expression on the left has the form $\gamma^{a_1}\tau^{b_1}\gamma^{a_2}\tau^{b_2}=(0,\ldots,0,a_2,0,\ldots,0a_1),\tau^{b_1+b_2}$. Since $b_1$ is not divisible by $n$, then $a_2$ is in the position which differs from $n$, so the first component is nonzero and the equality could not be true. 
\par\medskip

	If $p=3$, then $\gamma^{a_1}\tau^{b_1}\gamma^{a_2}\tau^{b_2}\gamma^{a_3}\tau^{b_3}=(0,\ldots,0,a_2,0,\ldots,0a_1),\tau^{b_1+b_2}\ast(0,\ldots,a_3),\tau^{b_3}$. Depending on $b_3$, the number $a_3$ either sums with $a_1$, or with $a_2$, either doesn’t sum with them, but in any case the first component will be nonzero, so the equality could not be true.
\par\medskip

	We found that $p\ge 4$. Easy to check that any relation of the length $p=4$ is the consequence of the relation $\gamma\tau^l\gamma\tau^{-l}=\tau^l\gamma\tau^{-l}\gamma$. So we consider the relations with the length $p\ge 5$ and will show, that one can reduce their length, thus proving that they are the consequence of the found relations.
\par\medskip 

	The relation $\gamma^{a_1}\tau^{b_1}\ldots\gamma^{a_p}\tau^{b_p}=e$ imposes some conditions on the numbers $b_1,\ldots,b_p$. Note, that in order to have true equality it is necessary to exist such $1 \le k \le n-1$, that the number $b_1+\ldots+b_k$ is divisible by $n$. If it is not true then the first component of the multiplier $\gamma^{a_1}\tau^{b_1}$ doesn’t sum with any other component, and the first component will be nonzero. So such $k\ge 1$ exists. Moreover, $k\ge 3$, otherwise if $k\le 2$ then one can reduce the length. We have the chain of the equalities $\gamma^{a_1}\tau^{b_1}\gamma^{a_2}\tau^{b_2}\gamma^{a_3}\tau^{b_3}\ldots\gamma^{a_k}\tau^{b_k}\gamma^{a_{k+1}}\tau^{b_{k+1}\ldots}=\\ \gamma^{a_1}\tau^{b_1}\gamma^{a_2}\tau^{-b_1}\tau^{b_1+b_2}\gamma^{a_3}\tau^{-b_1-b_2}\tau^{b_1+b_2+b_3}\ldots\tau^{b_1+\ldots+b_{k-1}}\gamma^{a_k}\tau^{-b_1-\ldots-b_{k-1}}\\ \tau^{b_1+\ldots+b_{k-1}+b_k}\gamma^{a_{k+1}}\tau^{b_{k+1}}\ldots\gamma^{a_p}\tau^{b_p}$. 
\par\medskip

	As $b_1+\ldots+b_k$ is divisible by $n$, then $\tau^{b_1+\ldots+b_k}=e$, so one can rewrite the expression in the following form $\\ \gamma^{a_1}\tau^{b_1}\gamma^{a_2}\tau^{-b_1}\tau^{b_1+b_2}\gamma^{a_3}\tau^{-b_1-b_2}\tau^{b_1+b_2+b_3}\ldots\tau^{b_1+\ldots+b_{k-1}}\gamma^{a_k}\tau^{-b_1-\ldots-b_{k-1}}\gamma^{a_{k+1}}\tau^{b_{k+1}}\ldots=\\ \gamma^{a_1+a_{k+1}}\tau^{b_1}\gamma^{a_2}\tau^{b_2}\ldots\gamma^{a_k}\tau^{-b_1-b_2-\ldots-b_{k-1}+b_{k+1}}\ldots\gamma^{a_p}\tau^{b_p}$. Here we shifted $\gamma^{a_{k+1}}$ to the left, using the relation $\gamma\tau^l\gamma\tau^{-l}=\tau^l\gamma\tau^{-l}\gamma$. The length of the resulting expression decreased at least by $1$, so it is the consequence of the found relations.
\par\medskip

	The last case, when the expression contains $\tau$, $\sigma$ and $\gamma$. The general form of a relation is $\gamma^{\alpha}\tau^{\beta}\sigma\ldots\gamma^{\delta}\ldots=e$. The main idea is that this expression is equivalent to the expression of the form $\sigma\tau^{u_1}\sigma\tau^{u_2}\ldots\sigma\tau^{u_s}\gamma^{a_1}\tau^{b_1}\gamma^{a_2}\tau^{b_2}\ldots\gamma^{a_p}\tau^{b_p}=e$, that is $\gamma$ and $\tau$ are gathered on the right, but $\sigma$ and $\tau$ are on the left. Indeed, if the expression has a fragment $\gamma^a\tau^b\sigma^c$, then for $0 \le \widehat b\le n-3$ one can rewrite this fragment in the next form $\gamma^a\tau^b\sigma=\tau^b\sigma\tau^{-b}\gamma^a\tau^{-b}$, where the remainder of the division of $b$ by $n$ is designated by $\widehat b$, i.e. $\gamma$ moves to the right, and $\sigma$ moves to the left. 
\par\medskip

	Now we need to understand the situation, when $\widehat b=n-2$ and $\widehat b=n-1$. Let $b=-1$, then the corresponding fragment is equivalent to $\gamma^a\tau^{-1}\sigma=\tau^{-1}\sigma\tau^2\gamma^a\tau^{-2}$, that is $\sigma$ moves to the left, and $\gamma$ moves to the left. If $b=-2$, then the corresponding fragment is equivalent to $\gamma^a\tau^{-2}\sigma=\tau^{-2}\sigma\tau\gamma\tau^{-1}$, so the situation is similar to the previous.
\par\medskip

	 Easy to see, that the condition $\sigma\tau^{u_1}\sigma\tau^{u_2}\ldots\sigma\tau^{u_s}\gamma^{a_1}\tau^{b_1}\gamma^{a_2}\tau^{b_2}\ldots\gamma^{a_p}\tau^{b_p}=e$ splits into the two following conditions: $\sigma\tau^{u_1}\sigma\tau^{u_2}\ldots\sigma\tau^{u_s}=\tau^{\alpha}$ and $\gamma^{a_1}\tau^{b_1}\gamma^{a_2}\tau^{b_2}\ldots\gamma^{a_p}\tau^{b_p}=\tau^{\beta}$. But as we have shown above these conditions both follow from the found relations. {\bfseries The assertion is proved.}
\par\medskip

	On the basis of these relations one can derive the group description $\mathbb Z^n\rtimes S_n$ via two generators, namely the next result holds:
\par\medskip

	{\bfseries Theorem 1.} For $n\ge 4$ the group $\mathbb Z^n\rtimes S_n$ has the following description:

$$
\left\{
\begin{array}{rcl}
b^2&=&e\\
(baba^{-1})^3&=&e\\
(ba^kba^{-k})^2&=&e,\ 2\le k\le n-2\\
ba^nba^{-n}&=&e\\
\end{array}
\right.
$$ 
\par\medskip

	{\bfseries Proof.} As the new generators we take $\\ a=(0,\ldots,0,1),(n12\ldots n-1)$ and $b=(0,\ldots,0),(213\ldots n)$, which are expressed though $\tau$, $\sigma$ and $\gamma$ by the next way:

$$
\left\{
\begin{array}{rcl}
a &=&\gamma\tau\\
b &=&\sigma.\\
\end{array}
\right.
$$
\par\medskip

       The relations between $\tau$, $\sigma$ and $\gamma$ are known from the previous assertion. Using them, we will obtain step by step the relations between $a$ and $b$. The relation $\sigma^2=e$ immediately implies the equaity $b^2=e$, so we found the first relation.
\par\medskip

	Rewrite the equality $\tau^m\sigma\tau^{-m}=\gamma\tau^m\sigma\tau^{-m}\gamma^{-1}=(\gamma\tau)\tau^{m-1}\sigma\tau^{-m+1}(\tau^{-1}\gamma^{-1})=a\tau^{m-1}\sigma\tau^{-m+1}a^{-1}$. The resulting equality implies $\tau^{-1}\sigma\tau=a^{-1}ba$ and $\tau^m\sigma\tau^{-m}=a^mba^{-m}$, $1\le m\le n-3$. With the help of this equality we can rewrite the relation $e=(\sigma\tau\sigma\tau^{-1})^3=(baba^{-1})^3$. So, we got the second relation.
\par\medskip

	Using the same method we obtain the third relation $e=(\sigma\tau^k\sigma\tau^{-k})^2=(ba^kba^{-k})^2$, but the range of changing for $k$ is slightly another, $2\le k\le n-3$. We need to get the missing equality $ba^{n-2}ba^{2-n}ba^{n-2}ba^{2-n}=e$. For this, we rewrite the equality $e=\sigma\tau^{n-2}\sigma\tau^{2-n}\sigma\tau^{n-2}\sigma\tau^{2-n}\\ =\sigma\tau\tau^{n-3}\sigma\tau^{3-n}\tau^{-1}\sigma\tau\tau^{n-3}\sigma\tau^{3-n}\tau^{-1}=baa^{n-3}ba^{3-n}a^{-1}baa^{n-3}ba^{3-n}a^{-1}=\\ ba^{n-2}ba^{2-n}ba^{n-2}ba^{2-n}$.
\par\medskip

	Let us now see what follows from $(\sigma\tau)^{n-1}=\tau^n$. One can rewrite it $\tau=\tau^{-1}(\sigma\tau)^{n-1}\tau^{2-n}=\tau^{-1}\sigma\tau\sigma\tau\sigma\tau^{-1}\tau^2\sigma\tau^{-2}\ldots\tau^{n-3}\sigma\tau^{3-n}=\\ a^{-1}bababa^{-1}a^2ba^{-2}\ldots a^{n-3}ba^{3-n}=a^{-1}(ba)^{n-2}ba^{3-n}$. Thus, the given equality implies the expression for $\tau$ via $a$ and $b$.
\par\medskip

	The equality $\gamma\tau^m\sigma\tau^{-m}=\tau^m\sigma\tau^{-m}\gamma$, $0\le m\le n-3$ holds automatically by substituting there the new generators and using the found relations, because we derived them by means of this equality.
\par\medskip

	It remains to ascertain what gives the last equality $\gamma\tau^{m}\gamma\tau^{-m}=\tau^{m}\gamma\tau^{-m}\gamma$. First let $m=1$, then $\gamma\tau\gamma\tau^{-1}=\tau\gamma\tau^{-1}\gamma$. Rewrite the equality $\tau\gamma\tau^{-1}=\gamma\tau\gamma\tau^{-1}\gamma^{-1}=a\gamma a^{-1}$. Substituting there the expression for $\tau$ and $\gamma$, we obtain $a^nb(a^{-1}b)^{n-2}=(ba)^{n-2}bab(a^{-1}b)^{n-2}a^{n-2}b(a^{-1}b)^{n-2}a$.
\par\medskip

	To simplify this equality we need the next lemma.
\par\medskip

	{\bfseries Lemma 4.} For the generators $a$ and $b$, taking into account the found relations, the relation $(a^{-1}b)^{n-k}a^{n-k}ba^{-1}=a^{-1}baba^{-2}b(a^{-1}b)^{n-k-2}a^{n-k-1}$ holds.
\par\medskip

	{\bfseries Proof.} Using the third relation, we can rewrite $$(a^{-1}b)^{n-k}a^{n-k}ba^{-1}=(a^{-1}b)^{n-k-1}a^{n-k-1}ba^{k-n}ba^{n-k-1}=(a^{-1}b)^{n-k-2}a^{n-k-2}ba^{1+k-n}ba^{-1}ba^{n-k-1}=$$ 
$$\ldots=a^{-1}baba^{-2}b(a^{-1}b)^{n-k-2}a^{n-k-1}.$$ {\bfseries The lemma is proved.}
\par\medskip

	Using the previous lemma for $4\le k\le n-2$ let us rewrite $\\ (ba)^{n-k}ba^{k-1}b(a^{-1}b)^{n-2}a^{n-k}b(a^{-1}b)^{n-k}a=\\ (ba)^{n-k}ba^{k-1}b(a^{-1}b)^{k-2}(a^{-1}b)^{n-k}a^{n-k}b(a^{-1}b)^{n-k}a=\\ (ba)^{n-k}ba^{k-1}b(a^{-1}b)^{k-2}a^{-1}baba^{-2}b(a^{-1}b)^{n-k-2}a^{n-k-1}b(a^{-1}b)^{n-k-1}a=\\ (ba)^{n-k}ba^{k-1}b(a^{-1}b)^{k-3}ba^{-2}baba^{-1}b(a^{-1}b)^{n-k}a^{n-k-1}b(a^{-1}b)^{n-k-1}a=\ldots=(ba)^{n-k}ba^{k-1}ba^{-1}ba^{2-k}ba^{k-3}b(a^{-1}b)^{n-4}a^{n-k-1}b(a^{-1}b)^{n-k-1}a=\\ (ba)^{n-k-1}ba^{k}ba^{1-k}ba^{k-2}ba^{2-k}ba^{k-3}b(a^{-1}b)^{n-4}a^{n-k-1}b(a^{-1}b)^{n-k-1}a=\\ (ba)^{n-k-1}ba^{k}b(a^{-1}b)^{n-2}a^{n-k-1}b(a^{-1}b)^{n-k-1}a$. Easy to check that this relation holds also for $k=2$ and $k=3$.
\par\medskip

	Using the equality above $n-3$ times, we obtain $\\ (ba)^{n-2}bab(a^{-1}b)^{n-2}a^{n-2}b(a^{-1}b)^{n-2}a=\\ baba^{n-2}b(a^{-1}b)^{n-2}aba^{-1}ba=baba^{n-2}b(a^{-1}b)^{n-4}a^{-2}bab=\\ baba^{n-2}b(a^{-1}b)^{n-5}a^{-3}ba^2ba^{-1}b=\ldots=baba^{n-2}ba^{2-n}ba^{n-3}b(a^{-1}b)^{n-4}=\\ ba^{n-1}b(a^{-1}b)^{n-3}$. As a result, the initial equality becomes $a^nb(a^{-1}b)^{n-2}=ba^{n-1}b(a^{-1}b)^{n-3}$. Reducing, we get $a^nba^{-1}=ba^{n-1}$, or $a^nb=ba^n$, what was required to get.
\par\medskip

       Now let $m\ge 2$, we have the equality $\gamma\tau^m\gamma\tau^{-m}\gamma^{-1}=\tau^m\gamma\tau^{-m}$. Rewrite it $\gamma\tau^m\gamma\tau^{-m}\gamma^{-1}=(\gamma\tau)\tau^{m-1}\gamma\tau^{-m+1}(\tau^{-1}\gamma^{-1})=a\tau^{m-1}\gamma\tau^{-m+1}a^{-1}=aa^{m-1}\gamma a^{-m+1}a=a^m\gamma a^{-m}$. We need to check that for $m\ge 2$ the equality $\tau^m\gamma\tau^{-m}=a^m\gamma a^{-m}$ holds automatically by the relations between $a$ and $b$ which we found above. For $m=1$ this condition holds. Let this condition holds for $m$, then we must prove, that it holds for $m+1$. We have $\tau^{m+1}\gamma\tau^{-m-1}=\tau\tau^m\gamma\tau^{-m}\tau^{-1}=\tau a^m\gamma a^{-m}\tau^{-m}=\tau a^{n+m-2}b(a^{-1}b)^{n-2}a^{1-m}\tau^{-1}=a^{-1}(ba)^{n-2}ba^{m+1}b(a^{-1}b)^{n-2}a^{n-m-2}b(a^{-1}b)^{n-2}a=a^{-1}(ba)^{m}(ba)^{n-m-2}ba^{m+1}b(a^{-1}b)^{n-2}a^{n-m-2}b(a^{-1}b)^{n-m-2}aa^{-2}b(a^{-1}b)^{m-1}a=a^{-1}(ba)^mba^{n-1}b(a^{-1}b)^{n-3}a^{-2}b(a^{-1}b)^{m-1}a=a^{n-1}(ba)^{m-1}bab(a^{-1}b)^{n-2}a^{-2}b(a^{-1}b)^{m-1}a=a^{n-1}baba^{m-1}b(a^{-1}b)^{n-2}a^{-m}ba^{-1}ba$. The last equality is derived similarly to the calculations above.
\par\medskip

	We need to establish the equality $a^{n-1}baba^{m-1}b(a^{-1}b)^{n-2}a^{-m}ba^{-1}ba=a^{n+m-1}b(a^{-1}b)^{n-2}a^{-m}$, which is equivalent to $baba^{m-1}b(a^{-1}b)^{n-2}a^{-m}ba^{-1}=a^{m}b(a^{-1}b)^{n-2}a^{-m-1}b$. Rewrite the left, we obtain $\\ baba^{m-1}b(a^{-1}b)^{n-2}a^{-m}ba^{-1}=baba^{m-1}b(a^{-1}b)^{m-2}(a^{-1}b)^{n-m}a^{-m}ba^{-1}\\ =baba^{m-1}b(a^{-1}b)^{m-2}a^{-1}baba^{-2}b(a^{-1}b)^{n-m-2}a^{-m-1}=baba^{m-1}b(a^{-1}b)^{m-3}a^{-2}bab(a^{-1}b)^{n-m}a^{-m-1}=baba^{m-1}ba^{-1}ba^{2-m}ba^{m-3}b(a^{-1}b)^{n-4}a^{-m-1}=ba^{m}ba^{1-m}ba^{m-2}ba^{2-m}ba^{m-3}b(a^{-1}b)^{n-4}a^{-m-1}=ba^{m}b(a^{-1}b)^{n-2}a^{-m-1}$. Now the equality, the truth of which we need to establish, takes the form $ba^{m}b(a^{-1}b)^{n-2}a^{-m-1}=a^{m}b(a^{-1}b)^{n-2}a^{-m-1}b$, or $ba^{m}b(a^{-1}b)^{n-2}a^{-m-2}=a^{m}b(a^{-1}b)^{n-2}a^{-m-1}ba^{-1}$. Rewrite the right $\\ a^{m}b(a^{-1}b)^{n-2}a^{-m-1}ba^{-1}=a^{m}(a^{-1}b)^{m-1}(a^{-1}b)^{n-m-1}a^{-m-1}ba^{-1}=\\ a^{m}b(a^{-1}b)^{m-1}a^{-1}baba^{-2}b(a^{-1}b)^{n-m-3}a^{-m-2}=a^{m}b(a^{-1}b)^{m-2}a^{-2}bab(a^{-1}b)^{n-m-1}a^{-m-2}=a^{m}ba^{-1}ba^{1-m}ba^{m-2}b(a^{-1}b)^{n-4}a^{-m-2}=a^mba^{-m}ba^{m-1}b(a^{-1}b)^{n-3}a^{-m-2}=ba^{m}b(a^{-1}b)^{n-2}a^{-m-2}$. Thus, we have established the truth of the equality. {\bfseries The theorem is proved.}
\par\medskip

	{\bfseries Notation 3.} The theorem 1 does not include two cases, namely, when $n=2$ and $n=3$. It holds for $n\ge 4$ because the standard description of the group $S_n$ via transpositions holds for $n\ge 4$. For the cases $n=2$ and $n=3$ some modifications are necessary.
\par\medskip

       For $n=3$ we have the description $\mathbb Z^3\rtimes S_3$ 
$$
\left\{
\begin{array}{rcl}
b^2&=&e\\
(baba^{-1})^3&=&e\\
ba^3ba^{-3}&=&e\\
\end{array}
\right.
$$
\par\medskip

	For $n=2$ we have the description $\mathbb Z^2\rtimes S_2$
$$
\left\{
\begin{array}{rcl}
b^2&=&e\\
ba^2ba^{-2}&=&e\\
\end{array}
\right.
$$	 
\par\medskip

{\bfseries\Large 3. The geometric realization of the groups}

\par\medskip
\par\medskip

{\bfseries Assertion 4.} In the $n$-dimensional space the points of $\mathfrak A_n$ represent the vertices of the prisms with the $(n-1)$-dimensional permutohedron base, forming the tessellation of the space.
\par\medskip

{\bfseries Proof.} In assertion $2$ we translated the points of $\mathfrak A_n$ by vector $s=(0,n-1,\ldots,2,1)$ and obtained all such points, the coordinates of which are distinct modulo $n$, and the geometric picture does not change under translating the points. The data points are of the form $(nm_1+l_1,\ldots,nm_n+l_n)$, where $l_1,\ldots,l_n$ are the numbers $1,2,\ldots,n$ in some order. 
\par\medskip

	The $n-1$-dimensional permutohedron is the convex hull of $n!$ points, which are obtained from the point $(1,2,\ldots,n)$ by permuting its coordinates. According to theorem $2$ from [2], the vectors $e_1=(-(n-1),1,\ldots,1,1),\ e_2=(1,-(n-1),\ldots,1,1),\ldots,\ e_n=(1,1,\ldots,1,-(n-1))$ translate the permutohedron in the hyperplane $x_1+\ldots+x_n=\frac{n(n+1)}{2}$, such that its parallel copies tessellate this plane. These vectors form the lattice, with respect to which the permutohedron is the polytope of Voronoi, so the parallel copies of permutohedron either don’t intersect or intersect in a common face. The vector $e_n$ is the leaner combination of the vectors $e_1,\ldots,e_{n-1}$ since $e_1+\ldots+e_n=0$. Let us take the system $e_1,e_2,\ldots,e_{n-1},a$, where $a=(1,1,\ldots,1)$. If we translate the initial permutohedron by all the integer leaner combinations of these vectors, then we obviously obtain the prisms with the permutohedron base, such that they form the tessellation of the $n$-dimensional space. Let us prove, that the set of the vertices of these prisms (denote this set by $\mathfrak B$) coincides with the set $s+\mathfrak A_n$.
\par\medskip

	Any point of the set $\mathfrak B$ has the general form $(-(n-1)t_1+t_2\ldots+t_n,\ldots,t_1+\ldots-(n-1)t_{n-1}+t_n,t_1+\ldots+t_n)+(u_1,\ldots,u_n)$, where $u_1,\ldots,u_n$ are the numbers $1,2,\ldots,n$ in some order. One can see, that all the coordinates of such points are equal modulo $n$, so $\mathfrak B\subset\mathfrak A_n$. Let us prove the reverse inclusion. To do it, we have to solve the system of the linear equations $(-(n-1)t_1+t_2\ldots+t_n,\ldots,t_1+\ldots-(n-1)t_{n-1}+t_n,t_1+\ldots+t_n)+(u_1,\ldots,u_n)= (nm_1,\ldots,nm_n)+(l_1,\ldots,l_n)$ with respect to $t_1,\ldots,t_n,u_1,\ldots,u_n$. Let $u_1=l_1,\ldots,u_n=l_n$, then the system takes the form 
$$
\begin{pmatrix}
nm_1\\
nm_2\\
\ldots\\
nm_n
\end{pmatrix}
=C\begin{pmatrix}
t_1\\
t_2\\
\ldots\\
t_n
\end{pmatrix}
,$$
$$
C=\begin{pmatrix}
-(n-1) & 1 & 1 & \ldots & 1\\
1 & -(n-1) & 1 & \ldots & 1\\
\vdots & \vdots & \ddots & \vdots & \vdots\\
1 & 1 & \ldots & -(n-1) & 1\\
1 & 1 & 1 & \ldots &1
\end{pmatrix}
.$$
\par\medskip

The inverse matrix can be calculated explicitly

$$
C^{-1}=\begin{pmatrix}
-\frac{1}{n} & 0 & 0 & \ldots & \frac{1}{n}\\
0 & -\frac{1}{n} & 0 & \ldots & \frac{1}{n}\\
\vdots & \vdots & \ddots & \vdots & \vdots\\
0 & 0 & \ldots & -\frac{1}{n} & \frac{1}{n}\\
\frac{1}{n} & \frac{1}{n} & \frac{1}{n} & \ldots & \frac{1}{n}
\end{pmatrix}
.$$
\par\medskip

	Obviously, the system has the integral solution, so $\mathfrak A_n=\mathfrak B$. {\bfseries The assertion is proved.}
\par\medskip

	Now suppose that the action of the group $\mathbb Z$ on $V$ is defined by any permutation $\tau$ and let denote the group of the invertible elements by $\mathfrak T_n^\tau$. Every permutation is the product of the disjoint cycles in the unique way up to a permutation of the factors, so $\tau=\tau_1\ldots\tau_k$, where $\tau_i$ are cycles, $i=1,\ldots,k$. This representation also defines the partition of the set of $n$ elements in $k$ disjoint subsets. Denote this partition by $\sigma$. Clear, that the vector $(a_1\ldots,a_n)$ belongs to the group $\mathfrak T_n^\tau$ if and only if the elements of the partition $\sigma$ of this vector belong to the $\mathfrak A_{|\tau_i|}$, $i=1,\ldots,k$, where $|\tau_j|$ is the length of the cycle $\tau_j$, and the elements of the partition change independently.
\par\medskip

{\bfseries Theorem 2.} In the $n$-dimensional space the points of the set $\mathfrak T_n^\tau$ are the vertices of the polytopes $\prod_{|\tau_1|-1}\times\ldots\times\prod_{|\tau_k|-1}\times I^k$, that form the tessellation of the space. For $\prod_i$ denotes $i$-dimensional permutohedron. Moreover, $\mathfrak T_n^\tau\cong(\mathbb Z^{|\tau_1|}\rtimes S_{|\tau_1|})\times\ldots\times(\mathbb Z^{|\tau_k|}\rtimes S_{|\tau_k|})$.
\par\medskip

{\bfseries Proof.} According to the assertion $4$, the groups $\mathfrak A_{|\tau_i|}$ in $\mathbb Z^{|\tau_i|}$ represent the vertices of the tessellation with $\prod_{|\tau_i|-1}\times I$, that is the prisms with the permutohedron base. Then every vector $\mathfrak T_n^\tau$ in $\mathbb Z^n=\mathbb Z^{|\tau_1|}\oplus\ldots\oplus Z^{|\tau_k|}$ is obtained from the elements in $\mathfrak A_{|\tau_i|}\subset\mathbb Z^{|\tau_i|}$ taken one by one from each group and assembled with the partition $\sigma$. Therefore, for every polytope $P_i$ of the partition $\mathbb Z^{|\tau_i|}$, $i=1,\ldots ,k$ we obtain the product $P_1\times\ldots\times P_k$ in $\mathbb Z^n$. 
\par\medskip

	All the polytopes in $\mathbb Z^{|\tau_i|}$ are obtained from a fixed one by the parallel shift, so the resulting product differs from a fixed by the parallel shift. It is also clear that the obtained from the direct product polytopes cover the whole space $\mathbb Z^n$ with no overdubs, and polytopes either disjoint or intersect in a common face and form the required tessellation. 
\par\medskip

	The isomorphism is obvious from the arguments above. {\bfseries The theorem is proved.}
\par\medskip

{\bfseries References}
\par\medskip

[1] V. M. Buchstaber, “Semigroups of maps into groups, operator doubles, and complex cobordisms”, Topics in topology and mathematical physics, Amer. Math. Soc. Transl. Ser. 2, 170, Amer. Math. Soc., Providence, RI, 1995, 9–31.
\par\medskip

[2] Garber, A. I., Poyarkov, A.P. (2006), “On permutohedra”, Vestnik MGU, ser. 1, 2006, N 2, pp. 3-8. 

\end{document}